\documentclass[12pt]{article}

\usepackage{amsmath}
\usepackage{amsfonts}
\usepackage[T1]{fontenc} 
\usepackage{a4wide}

\newtheorem{theo}{Theorem}[section]
\newtheorem{prop}{Proposition}[section]
\newtheorem{corol}{Corollary}[section]
\newtheorem{df}{Definition}[section]
\newtheorem{lemme}{Lemma}[section]

\newcommand{\CQFD}{\hfill $\square$}

\newcommand{\ind}{\mathbf{1}}

\newcommand{\cC}{{\cal C}}
\newcommand{\cD}{{\cal D}}
\newcommand{\cE}{{\cal E}}
\newcommand{\cF}{{\cal F}}

\newcommand{\cS}{{\cal S}}

\newcommand{\bR}{{\bf R}}

\newfont{\msbm}{msbm10 scaled\magstep1}
\newfont{\msbms}{msbm7 scaled\magstep1} 

\newcommand{\bbD}{\mbox{$\mbox{\msbm D}$}}
\newcommand{\bbE}{\mbox{$\mbox{\msbm E}$}}

\newcommand{\bbP}{\mbox{$\mbox{\msbm P}$}}

\newcommand{\bbR}{\mbox{$\mbox{\msbm R}$}}

\newcommand{\bbZ}{\mbox{$\mbox{\msbm Z}$}}

\newcommand{\bbsR}{\mbox{$\mbox{\msbms R}$}}

\newcommand{\bbsZ}{\mbox{$\mbox{\msbms Z}$}}

\def\E{\mathop{\hbox{\rm I\kern-0.20em E}}\nolimits}

\def\d{{\rm d}}

\def\og{\leavevmode\raise.3ex
     \hbox{$\scriptscriptstyle\langle\!\langle$~}}
\def\fg{\leavevmode\raise.3ex
     \hbox{~$\!\scriptscriptstyle\,\rangle\!\rangle$}~}
     
\begin{document}

\title{Discrete approximation of stable white noise - Application to  spatial  linear filtering.}
\author{ Clément Dombry \protect\hspace{1cm}}
\maketitle

\abstract{Motivated by the simulation of stable random fields, we consider the issue of discrete approximations of independently scattered stable noise. Two approaches are proposed: grid approximations available when the underlying space is $\bbR^d$ and shot noise approximations available on more general spaces. Limit theorems stating the convergence of discrete random noises to stable white noise are proved. 
These results are then applied to study moving average spatial random fields with heavy-tailed innovations and related limit theorems. A second application deals with discrete approximation for Brownian Lévy motion on the sphere or on the euclidean space.

}
{\bf Key words:} stable white noise, stable fractional noise, linear filtering. 

\section{Motivations}

Stable integration is a basic tool in the theory of stable random process. Indeed, if $W$ denotes a stable random noise on the measurable space $(E,\cE)$ with control measure $m$ and $(f_t)_{t\in T}$ a kernel such that $f_t\in L^\alpha(E,\cE,m)$ for any $t\in T$, then the random process defined by
$$X_t=\int_E f_t(x) W(dx), \quad t\in T $$
is a stable process, and the path properties of the process can be deduced from the properties of the kernel (see for instance \cite{ST94} chapter 10).
For example, the first examples of  self-similar stationary increments (SSSI) stable process were constructed in such a way: using a moving average kernel, Taqqu \& Wolpert \cite{TW83} and Maejima \cite{Ma83} constructed the linear fractional stable motion;  the harmonisable fractional stable motion was proposed by Cambanis \& Maejima \cite{CM89} as the stable integral of the harmonisable kernel. More recently, Cohen and Samorodnitsky \cite{CS} propose a new class of SSSI stable process defined as the stable integral of a random kernel, the random kernel being the local time of a fractional Brownian motion. Generalizing the notion of selfsimilarity for random field in higher dimension, operator scaling stable random fields were defined and constructed by Biermé \& Scheffler \cite{BMS} as stable integrals of high-dimensional kernel with suitable scaling properties. \\
\ \\
For the purpose of simulation of these processes, we need a general theory of discrete approximation of stable random measures. 
An approach based on Lepage's series also called shot noise series was developped: generalized shot noise series were introduced for simulation in \cite{Ro87}, further developments were done in \cite{Ro90} and \cite{Ro01} and a
general framework was developed in \cite{CLL08}. In this paper, we propose a new approach for the discrete approximation of stable random noise. We propose and developp two alternatives: in the case when $W$ is a stable white noise on $\bbR^d$, a grid approximation can be used (see Theorem \ref{theo:whitecase}); for more general spaces, an approximation of the stable white noise by a Poisson random measure is discussed (see Theorem \ref{theo:poisson}).  
\ \\
We apply our results on grid approximation of stable random noise on $\bbR^d$ to the study of stable noise obtained by linear filtering from i.i.d. random fields in the domain of attraction of stable distributions. This extends the results of Kokoszka \& Taqqu \cite{KT94,KT94b,KT95,KT96} for one dimensional sequences such as ARMA or FARIMA with stable inovations to spatial random fields. We prove a limit theorem for the suitably rescaled random noise (see Theorem \ref{theo:fraccase}). Two cases occur: if the coefficients of the random filter decrease quickly, the limit random field is stable white noise and the dependence vanishes in the limit; if the coefficients decrease slowly like a power function, the limit random field is fractional stable random noise with long range dependence.\\
\ \\
Another application of the Poisson approximation for stable white noise is given in the framework of stable Lévy motion on the sphere $\cS^d$ or on the euclidean space $\bbR^d$. These processes mimick the simple covariance structure of standard Brownian
motion on $\bbR$ to more general metric spaces (see \cite{KTU81}). As a direct application of our results, we give a Donsker's type Theorem for stable Lévy motion on the sphere or on the euclidean space.\\
\ \\
The paper is organised as follows. In section 2, we remind the reader of general results on random noise, stable white noise, Poisson random measures and their convergence. In section 3, the convergence of discrete random measure to stable white noise is proven for two different models: grid approximations and Poisson approximation. In section 4, two applications are exposed: linear filtering of i.i.d. random field and Donsker's type theorem for Lévy stable motion on the sphere or on the euclidean space. Section 5 is devoted to the proof of our results. Technical results on convergence of deterministic functions are gathered in an appendix.

\section{Stable random noises and Poisson random measures}
A random noise on $(E,\cE)$ is a generalized random field $(W[f])_{f\in\cF}$ indexed by a linear subspace $\cF$ of the space of real-valued measurable functions on $E$, verifying the linearity property: for all $a_1,a_2\in\bbR$ and $f_1,f_2\in\cF$.
\begin{equation}\label{eq:lin}
W[a_1f_1+a_2f_2]=a_1W[f_1]+a_2W[f_2] \ \ {\rm a.s.}
\end{equation}
Using the Cramer-Wold's device, this implies that the finite-dimensional distributions of the process $(W[f])_{f\in\cF}$ are determined by 
its one-dimensional distributions. The set of functions $\cF$ is  thought as the set of integrands. To emphasize the analogy with usual integration, we use the notation
$$W[f]=\int_E f(x)W(dx).$$
If $\ind_A\in \cF$ for some measurable set $A\subseteq E$, the value $W(A)=W[\ind_A]$ is thought as the measure  of $A$. However, each realisation of $W$ does not define an usual signed measure because it needs not to have finite variations. That's why some author (including myself) prefer the term random noise to random measure, this last term being reserved for measure-valued random variables.

We recall in this section the definition and properties of stable white noises (including Gaussian white noises) and Poisson random measures.

\subsection{Independently scattered $\alpha$-stable random noises}
Let $(E,\cE,m)$ be a measured space with  $m$  a $\sigma$-finite measure. For $\alpha>0$, we denote by $L^\alpha=L^\alpha(E,\cE,m)$ the space of measurable functions $f:E\to\bbR$ such that $\int_E|f|^\alpha dm<\infty$. For $f\in L^\alpha$, let $|\!|f|\!|_{L^\alpha}=\left(\int_E |f|^\alpha dm\right)^{1/\alpha}.$
If $\alpha\leq 1$, this defines a norm and $L^\alpha$ is a Banach case. This is no longer the case if $0<\alpha<1$. 

The stable distribution of index $\alpha\in (0,2]$ and parameters $\sigma \geq 0$ (scale), $\nu\in [-1,1]$ (skewness) is denoted by $\cS_\alpha(\sigma,\nu)$. For the sake of simplicity, we will always assume that $\nu=0$ if $\alpha=1$. The distribution $\cS_\alpha(\sigma,\nu)$ is defined by its Fourier transform
\begin{equation}\label{eq:fcstable1}
 \bar\lambda(\theta)=\exp\left(-\sigma^\alpha|\theta|^\alpha\left(1-i \nu \varepsilon(\theta)\tan\frac{\pi\alpha}{2}\right)\right), \quad \theta\in\bbR
\end{equation}
where $\varepsilon(a)=+1$ if $a>0$, $\varepsilon(a)=-1$ if $a< 0$ and $\varepsilon(0)=0$. 

The $\alpha$-stable random noise $W_\alpha$ on $(E,\cE)$ with control measure $m$ and  skewness function $v:(E,\cE)\to [-1,1]$ is the stable random process  defined on some probability space $(\Omega,\cF,\bbP)$  and indexed by $L^\alpha$, verifying the linearity property \ref{eq:lin} and with one dimensional marginal distributions given by
$$ W_\alpha[f]\sim S_\alpha(\sigma_f,\nu_f), \quad f\in L^\alpha$$
 where  
\begin{equation}\label{eq:param}
 \sigma_f=\left(\int_E |f|^\alpha dm\right)^{1/\alpha}\quad \mbox{and}\quad  \sigma_f^\alpha \nu_f=\int_E \nu \varepsilon(f)|f|^{\alpha}dm.
\end{equation}
Evaluating the process $W_\alpha$ at functions $(f_i)_{1\leq i\leq k}$ with pairwise disjoint supports yields independent random variables $(W_\alpha[f_i])_{1\leq i\leq k}$: we say that the random noise  $M_\alpha$ is independently scattered or white to qualify this property.

In the Gaussian case $\alpha=2$, the parameter $\nu$ is irrelevant and we retrieve the usual Wiener integral, which is an isometry from $L^2(E)$ onto some Gaussian subspace of $L^2(\Omega,\cF,\bbP)$. In the case $1< \alpha < 2$, the stable integral induces an isometry from $L^\alpha(E)$ onto some subspace of stable random variables on $(\Omega,\cF,\bbP)$ equipped with the covariation norm (see \cite{ST94}).

\subsection{Poisson random measures}
Poisson random measures are measure-valued random variables and they induce a random noise on their sets of integrands. In the following presentation, we focus on the point of view of random noises, a general reference is \cite{K86}.

Let $(E,\cE,n)$ be a measured space with $n$ a $\sigma$-finite measure. We denote by $\cE_0$ the class of measurable sets with finite $n$-measure. A random measure $N$ on $(E,\cE)$ is called a Poisson random measure with intensity $n$ if and only if for any finite collection $(A_i)_{1\leq i\leq n}$ of pairwise disjoint sets in $\cE_0$, the random variables $(N(A_i))_{1\leq i\leq n}$ are independent with $N(A_i)$ following a Poisson distribution with mean $n(A_i)$. 
The ({\it a priori} random) space $L^1(E,\cE,N)$ of integrable functions with respect to $N$ is almost surely equal to 
$$\cF_N=\left\{f:E\to \bbR \ |\ |f|\wedge 1\in L^1(E,n)\right\}$$
and $N$ induces a random noise on $\cF_N$. We note  $N[f]=\int_E fdN$ for $f\in\cF_N$.
This random noise is characterized by the linearity property $(\ref{eq:lin})$ and the one dimensional marginals given by their characteristic functions
$$
\bbE\left[\exp(i\theta N[f]) \right]=\exp \left(\int_E  (e^{i\theta f}-1) dn\right), \quad f\in \cF_N,\quad \theta\in\bbR.
$$

\subsection{Convergence of random measures}
In this paper, we are mainly concerned with convergence of random noises $(\mu_h)_{h>0}$ to stable random noise $W$, either independently scattered (Theorem \ref{theo:whitecase}) or fractional (Theorem \ref{theo:fraccase}). We shall now give a precise definition of convergence of random noises.

For $h>0$, let $\mu_h$ be a random noise with integrand space $\cF_{\mu_h}$ and let $W$ be some random noise with integrand space $\cF_W$. Let furthermore $\cF$ be a subspace included in all the  $\cF_{\mu_h}$'s, $h>0$ and in $\cF_W$. 

\begin{df}
We say that $\mu_h$ converges weakly to $W$ on $\cF$ as $h\to 0$ and write
$$
\mu_h \overset{\cF}\Longrightarrow W \ \ {\rm as}\ \ h\to 0,
$$
if the finite dimensional marginal distributions of $(\mu_h[f])_{f\in\cF}$ weakly converge to those of $(W[f])_{f\in\cF}$.
\end{df}

Using the linearity property of the random noises and the Cramer-Wold's device, we easily see that convergence of finite dimensional distributions is equivalent to  convergence of  one dimensional distributions, which is in turn equivalent to  convergence of characteristic functions. Although  straightforward, this will be of constant use and is stated for future reference in the following proposition:

\begin{prop}\label{prop:critere}
The following statements are equivalent, where convergence is meant as $h\to 0$:
\begin{enumerate}
\item $\mu_h$ converges weakly to $W$ on $\cF$,
\item for any $f\in\cF$, $\mu_h[f]$ converges weakly to $W[h]$,
\item for any $f\in\cF$ and $\theta\in\bbR$, $\bbE\left[\exp(i\theta\mu_h[f])\right]$ converge to $\bbE\left[\exp(i\theta W[f])\right]$.
\end{enumerate}
\end{prop}

To illustrate this definition, we give a first example of weak convergence of random noises.\\

{\bf Example :}
Let $(E,d)$ be some metric space and $\cE$ the borelian $\sigma$-algebra. Let $m$ and $(m_h)_{h>0}$ be finite measures on $E$ and $\nu$ and $\nu_h$ be measurable functions $E\to [-1,1]$. Define $W$ (resp. $W_h$) as the $\alpha$-stable random measure with control measure $m$ (resp. $m_h$) and skewness function $\nu$ (resp. $\nu_h$). Let $\cC_b(E)$ denotes  the space of bounded continuous functions on $E$. Note that from the assumption that $m$ is finite, $\cC_b(E)\subset \cF_W=L^\alpha(E,\cE,m)$, and  similarly $\cC_b(E)\subset\cF_{W_h}=L^\alpha(E,\cE,m_h)$ for all $h>0$.

Then the following two statements are equivalent, where convergences are meant as $h\to 0$:
\begin{enumerate}
\item the random noises $W_h$ converge weakly to $W$ on $\cC_b(E)$,
\item the control measures $m_h$ converge weakly to $m$  and the skewness measures $\nu_h dm_h$ converge weakly to $\nu dm$.
\end{enumerate}
The proof is direct once we recall that the $\alpha$-stable distributions $S_\alpha(\sigma_h,\nu_h)$ converge to $S_\alpha(\sigma,\nu)$ (with $\sigma>0$)if and only if $(\sigma_h,\nu_h)\to (\sigma,\nu)$.

\section{Discrete approximations of white stable noises}
We introduce two methods for approximating an independently scattered stable random noise on $(E,\cE)$. The first one, available when the underlying space is $E=\bbR^d$ and the control measure is the Lebesgue measure, relies on a grid approximation $h\bbZ^d\subseteq \bbR^d$ when the span $h$ of the grid goes to zero. The second one, avalaible on a general space $E$, relies on the approximation of the stable random noise by  suitable Poisson random measures or shot noise, when the intensity $\lambda$ goes to infinity.

\subsection{Grid approximation of white stable noise on $\bbR^d$}\label{sec:grid}
In this section, $E=\bbR^d$ and $W_\alpha$ denotes the independently scattered $\alpha$-stable random noise on $\bbR^d$ with Lebesgue intensity and constant skewness $\nu$. We propose a discrete approximation of $W_\alpha$ based on the grid approximation $h\bbZ^d\subseteq \bbR^d$ with span $h>0$.

The construction uses  a family $\xi=\{\xi_k,k\in\bbsZ^d\}$  of real random variables satisfying the following assumption:
\begin{center}
${\bf H}_\alpha$ \ \ \ the $\xi_k$'s are i.i.d. random variables in the normal domain of attraction of the stable distribution $\cS_\alpha(\sigma,\nu)$.
\end{center}
We suppose furthermore $\nu=0$ if $\alpha=1$. Recall that the normal domain of attraction of the stable distribution $\cS_\alpha(\sigma,\nu,0)$ consists in the random variables $Y$ such that 
$$
n^{-1/\alpha}\sum_{i=1}^n Y_i \Longrightarrow \cS_\alpha(\sigma,\nu) {\rm\ as\ } n\to\infty 
$$ 
where the $Y_i$ are i.i.d. random variables distributed as $Y$.
If $\alpha=2$, ${\bf H}_2$ holds if and only if the $\xi_k$'s are i.i.d. random variables with $\bbE[\xi_k]=0$ and $\bbE[\xi_k^2]=2\sigma^2$. The stable distribution $S_2(\sigma,0)$ is then  the normal distribution with mean $0$ and variance $2\sigma^2$ (the skewness parameter $\nu$ is irrelevant in this case).
In the case $0<\alpha<2$, from the known characterization of the domain of attraction of stable distributions (see \cite{Feller}), the $\xi_k$'s must satisfy the following tail assumptions as $x\to +\infty$:
\begin{equation}\label{eq:heavy-tail}
\bbP(\xi_k\geq x)\sim px^{-\alpha} {\rm\ and\ } \bbP(\xi_k\leq -x)\sim qx^{-\alpha},
\end{equation}
where the constants $\sigma, \nu,p,q$ satisfy
\begin{eqnarray*}
\sigma^\alpha &=& (p+q)\int_0^\infty t^{-\alpha}\sin t\d t,\\
\sigma^\alpha \nu \tan(\pi\alpha/2) &=& (p-q)\int_0^\infty t^{-\alpha}(1-\cos t)\d t.
\end{eqnarray*}
Note that equation (\ref{eq:heavy-tail}) is equivalent to the fact that $\xi-\tau$ belongs to the normal domain of attraction of the $\alpha$-stable distribution $S_\alpha(\sigma,\nu)$ for some $\tau\in\bbR$, but we assume here furthermore that $\tau=0$.\\

We propose the discrete approximation of $W_\alpha$ obtained by replacing $W_\alpha(dx)$ by $\gamma_h\xi_k \ind_{h(k+I^d)}(x)dx$ on each cell $h(k+I^d)$, $k\in\bbZ^d$. Here $I^d=[0,1)^d$ and $\gamma_h=\sigma^{-1}h^{(\frac{1}{\alpha}-1)d}$. More formally, let $\mu_h=\mu_h(\xi)$ be the random measure on $\bR^d$ absolutely continuous with respect to Lebesgue measure with random density
\begin{equation}\label{eq:defmu}
\frac{d\mu_h}{dx}(x)=\gamma_h\sum_{k\in\bbsZ^d}\xi_k  \ind_{h(k+I^d)}(x). 
\end{equation}

The random signed measure $\mu_h$ defined by (\ref{eq:defmu}) induces a random noise on the integrands set 
\begin{equation}\label{eq:defFmu}
\cF_{\mu_h}=\left\{f\in L^1_{loc}(\bbR^d)\ ;\ \int_{\bbsR^d}f d\mu_h {\rm\ converges\ a.s.\ } \right\}.
\end{equation}
Note that in this definition, only the semi-convergence of the integral $\int_{\bbsR^d} fd\mu_h$ or equivalently of the series
$$
\gamma_h\sum_{x\in\bbsZ^d}  \xi_k \int_{h(k+I^d)}f(x)dx 
$$
is required almost surely. 

The following scaling relation is worth noting:
\begin{equation}\label{eq:scaling}
\mu_h[f(c.)]=c^{-\alpha^{-1}d}\mu_{ch}[f(.)], \quad c>0
\end{equation}
whenever these quantities are well-defined.

Our first result precise the integrands sets $\cF_{\mu_h}$. We need to introduce the space $\cD^\alpha(\bbR^d)$ of locally integrable function $f\in L^1_{loc}(\bbR^d)$ such that $f(x)=o(|x|^{-\eta})$ for some $\eta>\alpha^{-1}d$. 

\begin{prop}\label{prop:domain}
Let $0<\alpha\leq 2$ and suppose that $\xi$ satisfies assumption ${\bf H}_\alpha$.\\
For $h>0$, consider the random measure $\mu_h$ defined by equation (\ref{eq:defmu}). Then, 
\begin{enumerate}
\item If $1\leq \alpha \leq 2$, $L^\alpha(\bbR^d)\subset \cF_{\mu_h}$.
\item If $0<\alpha< 1$, $\cD^\alpha(\bbR^d) \subset \cF_{\mu_h}$.
\end{enumerate}
\end{prop}

We then consider weak convergence of the random  measures $(\mu_h)_{h>0}$. 
To unify the notation, we introduce  $\cF^\alpha=L^\alpha(\bbR^d)$ if $1\leq \alpha\leq 2$, and $\cF^\alpha=\cD^\alpha(\bbR^d)$ if $0<\alpha<1$. 

\begin{theo}\label{theo:whitecase}
Let $0<\alpha\leq 2$ and suppose that $\xi$ satisfies assumption ${\bf H}_\alpha$.\\
Then as $h\to 0$, $\mu_h$ converge weakly to $W_\alpha$ on $\cF^\alpha$. 
\end{theo}

\noindent
{\bf Remark:} In the case $1\leq \alpha\leq 2$, we obtain convergence on $L^\alpha(\bbR^d)$ the full natural set of integrands for $W_\alpha$. In the case $0<\alpha<1$, we have to restrict on the smaller subspace $\cD^\alpha(\bbR^d)$ because $L^\alpha$ is not included in $\cF_{\mu_h}$ (see the remark after Lemma \ref{lem:lem1.2} in the Appendix.)\\
\ \\
\noindent
{\bf Remark:} We prove in fact a slightly stronger result 
$$\mu_h[f_h]\Rightarrow W_\alpha[f] {\rm\ as\ } h\to 0 {\rm \ if\ } f_h\to f {\rm\ in\ } \cF_\alpha .$$
See Proposition \ref{prop:diagconv} for a precise statement including the definition of convergence in $\cD^\alpha(\bbR^d)$ if $0<\alpha<1$. This diagonal convergence can be used to prove the following interesting corollary. \\
Let $\cC(\bbR^d)$ denote the space of continuous functions on $\bbR^d$. 
\begin{corol}\label{cor1}
Let $0<\alpha\leq 2$ and suppose that $\xi$ satisfies assumption ${\bf H}_\alpha$.\\
Then the random signed measure $\tilde\mu_h$ on $\bbR^d$ defined by 
$$
\tilde{\mu}_h= \sigma^{-1}h^{\alpha^{-1}d}\sum_{k\in\bbsZ^d} \xi_k \delta_{hk}  
$$
converges weakly to $W_\alpha$ on $\cF^\alpha\cap \cC(\bbR^d)$  as $h\to 0$.
\end{corol}

These results can be used to propose simulation of $\alpha$-stable processes represented as stable integrals. Another approach developed in \cite{CLL08} is to use the representation of stable integral in Lepage's series. For the purpose of simulation, it is important to provide error bounds or speed of convergence for the method. We now consider the speed of the convergence of the random noise $\mu_h$ to the stable integrals $W_\alpha$ in Theorem \ref{theo:whitecase}. Better rates of convergence are expected if the random variables $\xi$ satisfy ${\bf H}_\alpha$ with distribution $\cS_\alpha(\sigma,\nu)$ rather than in the domain of attraction of the stable distribution. In this case, the discrete approximation $\mu_h[f]$ has the same distribution as
$$\bar\mu_h[f]= \int_{\bbsR^d}(\tilde\psi_h f)(x)W_\alpha(dx),$$
where $\tilde\psi_h f$ is the discrete approximation of $f$ obtained by replacing $f$ by its mean value on each cell $h(k+[0,1[^d), k\in\bbZ^d$ (see equation \ref{eq:approx} below). In this representation $\xi_k$ corresponds to $\sigma h^{-d/\alpha}W_\alpha[h(k+[0,1[^d)]\sim \cS_\alpha(\sigma,\nu)$.
As a consequence, the random variable can be represented on the same space and we can focus on $L^p$ convergence and give $L^p$ error bounds. Furthermore, during the simulation, the kernel $f$ must sometimes be replaced by some approximated kernel $f_M$. For example, when $f$ has infinite support, $f_M$ might be the restriction of $f$ on a bounded domain. Or $f_M$ might be a piecewise constant approximation of $f$ so that the mean value of $f_M$ on $h(k+[0,1[^d)$ reduces to $f_M(hk)$. 

\begin{prop}\label{prop:errorbound}
Let $0<\alpha\leq 2$ and $f_M\in\cF^\alpha$ be an approximation of the kernel $f\in\cF^\alpha$ and let
$$\bar{\mu}_h[f_M]= \int_{\bbsR^d}(\tilde \psi_h f_M)(x)W_\alpha(dx).$$
Then, for every $0<p<\alpha$,
$$\bbE\big[\big|\bar{\mu}_h[f_M]-W_\alpha[f]\big|^p \big]\leq C_{\alpha,p}|\!|\tilde\psi_hf_M-f|\!|_{L^\alpha}^{p}$$
and $C_{\alpha,p}$ is the absolute $p$-th moment of the distribution $\cS_\alpha(1,\nu)$.\\
In the case $\alpha=2$, the result holds for all $p>0$
\end{prop}
{\bf Remark:} Note that $\tilde\psi_h f_M\to f_M$ in $\cF^\alpha$ as $h\to 0$ and explicit error bounds can be given under suitable regularity condition on $f_M$.

\subsection{Shot noise approximation of white stable noise on a general space}
In this section, $(E,\cE, m)$ is a general measured space with $m$ a $\sigma$-finite measure and $W_\alpha$ denotes the independently scattered $\alpha$-stable random noise on $E$ with control measure $m$ and constant skewness $\nu$. We propose a discrete approximation of $W_\alpha$ based on shot noises, i.e. on Poisson random measures.

Let $G$ be some distribution on $\bbR$ belonging to the normal domain of attraction of the stable distribution $\cS_\alpha(\sigma,\nu)$.  For $\lambda>0$, let $N_\lambda(de,d\xi)$ be the Poisson random measure on $E\times\bbR$ with intensity measure $n_\lambda(de,d\xi)=\lambda m(de)G(d\xi)$. 
It induces a random noise $\mu_\lambda$ on $E$ by the heuristic formula
$$ \mu_\lambda(dx) =\gamma_\lambda \int_{E\times \bbsR} \xi\delta_e(dx) N_\lambda(de,d\xi),$$
where $\gamma_\lambda= \sigma^{-1}\lambda^{-\frac{1}{\alpha}}$.
More formally, let $\cF_{\mu_\lambda}$ be the set of functions $f:E\to \bbR$ such that 
\begin{equation}\label{eq:def_Fmulambda}
\int_{E\times\bbsR} |\xi f(e)|\wedge 1 \ n_\lambda(d\xi,de) <\infty .
\end{equation}
We consider the random noise $\mu_\lambda$ on $\cF_{\mu_\lambda}$ defined by
\begin{equation}\label{eq:def_mulambda}
\mu_\lambda[f] =\gamma_\lambda \int_{E\times \bbR} \xi f(e) N_\lambda(de,d\xi).
\end{equation}
Our result is the following:

\begin{theo}\label{theo:poisson}
Suppose $G$ belongs to the normal domain of attraction of the stable distribution $\cS_\alpha(\sigma,\nu)$.
Then, for any $\lambda>0$, $L^\alpha(E,\cE,m)\subseteq \cF_{\mu_\lambda}$ and as $\lambda\to \infty$, $\mu_\lambda$ converge weakly to $W_\alpha$ on $L^\alpha(E,\cE,m)$.
\end{theo}

{\bf Remark:} In the case when $m$ is finite and normalized to $m(E)=1$, we can replace the Poisson random measure $N_\lambda$ by the binomial random measure $\tilde N_n=\sigma^{-1}n^{-\frac{1}{\alpha}}\sum_{i=1}^n \delta_{(\xi_i,e_i)}$ where $(\xi_i,e_i)_{i\geq 1}$ is a sequence of i.i.d. random variables with distribution $G(d\xi)m(de)$. We then consider the random noise
$$\tilde\mu_n[f]=\sigma^{-1}n^{-\frac{1}{\alpha}}\sum_{i=1}^n \xi_if(e_i).$$
Similar results as those of Theorem \ref{theo:poisson} hold: for any $n\geq 1$, $L^\alpha(E,\cE,m)\subseteq \cF_{\tilde\mu_n}$ and  $\tilde\mu_n$ weakly converge  to $W_\alpha$ on $L^\alpha(E,\cE,m)$ as $n\to\infty$. The proof is very similar to that of 
Theorem \ref{theo:poisson} and will be omitted.

\section{Applications}
We propose two applications of the above results about convergence to stable random noise. The first application based on Theorem \ref{theo:whitecase} is an analysis of random noises arrising from linear filtering of i.i.d. random fields. Interessant fields appear in the limit that we call fractional stable random noise. The second application based on Theorem \ref{theo:poisson} proposes a Donsker's type theorem for stable Lévy field on the sphere or on the euclidean space.

\subsection{Spatial linear filtering}
In this section, we consider a random measure $\hat\mu_h$  of the form 
\begin{equation}\label{eq:defhatmu}
\frac{d\hat\mu_h}{dx}(x)=\hat\gamma_h\sum_{k\in\bbsZ^d}\hat\xi_k  \ind_{h(k+I^d)}(x), 
\end{equation}
where $\hat\gamma_h>0$ is a normalisation constant and the field $\hat\xi=(\hat\xi_k)_{k\in\bbsZ^d}$ is not any more i.i.d. but stationary with spatial dependence (possibly long range dependence). 
More precisely, we consider the case when the field $\hat\xi$ is obtained from a random field $\xi$ satisfying ${\bf H}_\alpha$ by linear filtering. The linear filter is given by $c=(c_k)_{k\in\bbsZ^d} \in\ell^\alpha$ and the random field $\hat\xi$ by
$$\hat\xi_k =(\xi\ast c)_k=\sum_{l\in\bbsZ^d}c_{k-l}\xi_l .$$
The corresponding assumption is 
\begin{center}
${\bf\hat H}_\alpha$ \ \ \ $\hat\xi=\xi\ast c$ with $\xi$ satisfying ${\bf H}_\alpha$ and $c\in\ell^\alpha$.
\end{center}
Here $\ell^\alpha=\ell^\alpha(\bbZ^d)$ is the space of sequence $u=(u_k)_{k\in\bbZ^d}$ such that $\sum_{k\in\bbsZ^d}|u_k|^\alpha<\infty$ and $\ast$ denotes the convolution product of sequences.
We will see that if $\xi$ satisfies ${\bf H}_\alpha$ and $c\in\ell^\alpha$ then $\xi\ast c$ is defined almost surely.

The integrands set associated to the random measure $\hat\mu_h$ is   
\begin{equation}\label{eq:defFhatmu}
\cF_{\hat\mu_h}=\left\{f\in L^1_{loc}(\bbR^d)\ ;\ \int_{\bbsR^d}f d\hat\mu_h {\rm\ converges\ a.s.\ } \right\}.
\end{equation}

We consider the convergence of the random measure $\hat\mu_h$ to some $\alpha$-stable noise. According to the behavior of $c$ at infinity, the limit field can be either a white stable noise or a fractional stable noise (to be defined below). In the former case, spatial dependence disappear in the limit and weak dependence holds, whereas in the latter case, spatial dependence persists in the limit and strong dependence holds. For the sake of simplicity, we focus on the case $1<\alpha\leq 2$. As in section \ref{sec:grid}, $W_\alpha$ denotes a stable white noise on $\bbR^d$ with Lebesgue control measure and constant skewness $\nu$. 

\begin{theo}\label{theo:fraccase}
Let $1<\alpha\leq 2$ and suppose that $\hat\xi$ satisfies assumption ${\bf \hat H}_{\alpha}$.\\
For $h>0$, consider $\hat\mu_h$ the random measure given by equation (\ref{eq:defhatmu}).
\begin{enumerate}
\item Suppose $c\in\ell^1$ and let $\hat\gamma_h=\sigma^{-1}h^{({\alpha}^{-1}-1)d} $ and $C=(\sum_{k\in\bbsZ^d} c_k)$.\\ 
Then for any $h>0$, $L^\alpha \subseteq \cF_{\hat\mu_h}$ and as $h\to 0$, $\mu_h$ converge weakly to $CW_{\alpha}$ on $L^\alpha$. 
\item Suppose that $c$ has regular variations at infinity of order $-\beta\in(-d,-\alpha^{-1}d)$ in the sense that
\begin{equation}\label{eq:asympt}
\lim_{t\to\infty} \sup_{|x|=1} \left|t^\beta c_{[tx]}-p(x)\right|=0,
\end{equation}
where $p:\bbR^d\setminus \{0\}\rightarrow \bbR$ is some locally integrable homogeneous function of order $-\beta$ (hyp ??). Let $\hat\gamma_h=\sigma^{-1}h^{\alpha^{-1}d-\beta}$.\\ 
Then for any $h>0$, $L^\alpha\cap L^1\subseteq \cF_{\hat\mu_h}$ and as $h\to 0$, $\mu_h$ converge weakly to $W_{\alpha,p}$ on $L^\alpha\cap L^1$, where $W_{\alpha,p}$ is the fractional stable random noise defined by
$$W_{\alpha,p}[f]=W_\alpha[f\ast \check p]\quad,\quad f\in L^1\cap L^\alpha  $$ 
\end{enumerate}
\end{theo}

We now give some interessant properties of the random noise $W_{\alpha,p}$.
\begin{prop}\label{prop:fracfield}
Let $1<\alpha\leq 2$ and $p$ be as in Theorem \ref{theo:fraccase}.\\
\begin{enumerate}
\item The random noise $W_{\alpha,p}$ on $L^\alpha\cap L^1$ is $\alpha$-stable, stationary and $\alpha^{-1}d-\beta+d$ self-similar.
\item In the case $\alpha=2$, $W_{2,p}$ is the centered Gaussian random noise with covariance function
$$ \bbE\left(W_{\alpha,p}[f]W_{\alpha,p}[g]\right)=\int_{\bbR^d\times\bbR^d}f(x)K(x,y)g(y)\d x\d y,\quad f,g\in L^1\cap L^2$$
where $K:\bbR^d\times\bbR^d\to\bbR$ is the kernel given by
$$ K(x,y)=\int_{\bbsR^d}p(x-z)p(y-z)\d z.$$
\end{enumerate}
\end{prop} 

{\bf Remark:} The notions of stationarity and selfsimilarity are invariance properties of the random noise under the transformation group of translations and dilatations respectively. The formal definition is given in equation ?? and ??.

{\bf Remark} Self-similarity can be stated using the notion of renormalisation groups. Let $\cF$ be a linear space of functions on $\bbR^d$ closed under dilatations, i.e. for any $h>0$ and $f\in\cF$, $f(h.)\in\cF$. Let $u\in\bbR$. For $h>0$ let $T_h^{(u)}$ be the transformation acting on a random noise $W$ on $\cF$ by
$$(T_h^{(u)}W)[f(.)]= h^{u}W[f(h.)]\quad,\quad f\in\cF.$$
The group relations $T_1^{(u)}=Id$ and $T_h^{(u)}\circ T_{h'}^{(u)}=T_{hh'}^{(u)}$ holds and the group of transformation $(T_h^{(u)})_{h>0}$ is called the normalisation group of index $u$. Our results have a nice interpretation in terms of normalisation group. Consider first Theorem \ref{theo:whitecase}. The scaling relation (\ref{eq:scaling}) implies that the random measure $\mu_h$ is equal to $T_{h}^{(\alpha^{-1}d)}\mu_1$ and the Theorem states that the random noise $T_h^{(\alpha^{-1}d)}\mu_1$ converges to $W_\alpha$ as $h\to 0$. Furthermore, $W_\alpha$ is a fixed point of the renormalisation group of index $\alpha^{-1}d$, i.e. $T_h^{(\alpha^{-1}d)}W_\alpha=W_\alpha$, meaning that $W_\alpha$ is $(\alpha^{-1}d)$-selfsimilar. In the same way, consider Theorem \ref{theo:fraccase} part $2$. Here $\hat\mu_h=T_h^{(\alpha^{-1}d-\beta+d)}\hat\mu_1$ and the Theorem states that $T_h^{(\alpha^{-1}d-\beta+d)}\hat\mu_1$ converges to $W_{\alpha,p}$ as $h\to 0$. Furthermore, $W_{\alpha,p}$ is a fixed point of the renormalisation group of index $\alpha^{-1}d-\beta+d$, i.e. $T_h^{(\alpha^{-1}d-\beta+d)}W_{\alpha,p}=W_{\alpha,p}$, meaning that $W_{\alpha,p}$ is $\alpha^{-1}d-\beta+d$-selfsimilar.

\subsection{Convergence to Lévy Brownian motion on the sphere or on the euclidean space}
The Brownian motion parametrized by a general metric space was introduced by P.L\'evy. It is defined as follows: let $O$ be a fixed point in a metric space $(M,d)$, a Brownian motion parametrized with the metric space $(M,d)$ and with origin $O$ is a centered Gaussian process $(B(m))_{m\in M}$ such that :\\
- $B(0)=0$ almost surely,\\
- $B(m)-B(m')$ has variance $d(m,m')$.\\
There does not always exist a Brownian motion on $(M,d)$ for an arbitrary metric space $(M,d)$ but constructions of such Lévy Brownian motion have been proposed in the case when $M$ is a sphere, a euclidean space or an hyperbolic space (see \cite{KTU81}). These constructions are based on Gaussian white noise and have been extended to stable white noise in order to get stable processes with interesting properties.\\

First consider the case when $M=S^q$ is the unit sphere in $\bbR^{q+1}$ and $d$ is the geodesic distance. Let $O$ be a fixed point on the sphere $S^q$, for example the northern pole and let $W_2$ be a Gaussian white noise on the sphere with control measure $ds$ the normalized uniform measure on the sphere $S^q$.  For $m\in S^q$, let $H_m$ be the hemisphere centered at $m$ defined by
$$H_m=\{m'\in S^q\ |\ d(m,m')\leq \frac{\pi}{2}\}.$$
Then the process $(B(m))_{m\in S^q}$ defined by
$$B(m)=\sqrt{\pi}W_2\left(H_O\Delta H_m\right)\quad ,\quad m\in S^q$$
is a Brownian motion on $(S^q,d)$. Here $\Delta$ denotes the symmetric difference. More generally, the Lévy stable motion on the sphere $S^q$ is defined by the same formula with the Gaussian white noise $W_2$ replaced by a symmetric $\alpha$-stable random measure $W_\alpha$ with control measure $ds$.

When the metric space $(M,d)$ is the euclidean space $\bbR^q$, a Brownian motion is constructed in the following way (Lévy-Chenstov construction). Let $W_2$ be a white Gaussian noise on $S^q\times\bbR_+^\star$ with control measure $ds\times dr$. A pair $(s,r)$ represents the hyperplane in $\bbR^q$ with equation $<\!x,s\!>=r$ and $S^q\times\bbR_+^\star$ is thought as the set of hyperplanes that do not contain the origin. For $m\in\bbR^q$ define $V_m$ as the set of hyperplanes that separate $O$ and $m$, i.e.
$$V_m=\{(s,r)\in S^q\times\bbR_+^\star\ ;\ 0<r<<\!s,m\!>\}.$$
Then the process $(B(m))_{m\in \bbR^q}$ defined by
$$B(m)=W_2\left(V_m\right)\quad ,\quad m\in \bbR^q$$
is a Brownian motion on $(R^q,d)$. Replacing the Gaussian white noise $W_2$ by the symmetric $\alpha$-stable random measure $W_\alpha$ with control measure $dsdr$, we obtain the symmetric $\alpha$-stable Lévy-Chenstov random field, which is a $\alpha^{-1}$-self-similar stationary increments process.

We apply our results on convergence of Poisson random measures to stable noises to obtain a Donsker's theorem for stable Lévy motion on the sphere $S^q$ or  on $\bbR^q$.

\begin{theo}
Suppose $G$ is in the normal domain of attraction of the distribution $\cS_\alpha(\sigma,0)$ for some $\sigma>0$. 
The following convergence hold in the sense of finite dimensional distributions as $\lambda\to\infty$.
\begin{enumerate}
 \item Let $N_\lambda(ds,d\xi)$ be a Poisson random measure on $S^q\times\bbR$ with intensity $\lambda ds G(d\xi)$. The random process $B_\lambda$ on $S^q$ defined by
$$B_\lambda(m)=\sqrt{\pi}\sigma^{-1}\lambda^{-1/\alpha}\int_{S^q\times\bbsR} \xi{\bf 1}_{H_O\Delta H_m}(s) N_\lambda(ds,d\xi) \quad ,\quad m\in S^q$$
weakly converges  to the symmetric $\alpha$-stable Lévy motion on the sphere $S^q$.
\item Let $N_\lambda(ds,dr,d\xi)$ be a Poisson random measure on $S^q\times\bbR_+^\star\times\bbR$ with intensity $\lambda ds dr G(d\xi)$. The random process $B_\lambda$ on $\bbR^q$ defined by
$$B_\lambda(m)=\sigma^{-1}\lambda^{-1/\alpha} \int_{S^q\times\bbsR_+^\star\times\bbsR} \xi{\bf 1}_{V_m}(s,r) N_\lambda(ds,dr,d\xi) ,\quad m\in \bbR^q$$
weakly converges  to the symmetric $\alpha$-stable Lévy motion on the euclidean space $R^q$.
\end{enumerate}
\end{theo}

{\bf Remark :} The random fields $B_\lambda$ can be easily simulated, which is not clear for the limit field $B$. For example in the case of Lévy motion on the sphere $S^q$, simulating the Poisson integral can be made as follows. Draw $T$ according to a Poisson distribution with mean $\lambda$, conditionally to $T=t$, draw $(s_i,\xi_i)_{1\leq i\leq t}$ identically distributed with distribution $dsG(d\xi)$ and let 
$$ B_\lambda(m)=\sqrt{\pi}\sigma^{-1}\lambda^{-1/\alpha}\sum_{i=1}^t \xi_i{\bf 1}_{H_O\Delta H_m}(s_i).$$
Note that the quantity $x\in H_m$ is easily calculated since $x\in H_m$ if and only if $<\!x,m\!>\geq 0$. 
In the case of Lévy motion on the sphere $\bbR^q$, the method is just the same once we observe that simulating $B_\lambda(m), m\in D$ in a bounded domain $D$ requires the knowledge of the Poisson random measure on a bounded domain of $S^q\times\bbR^\star_+$ and hence only of a finite number of random points.


\section{Proofs}
\subsection{Preliminaries on stable distributions}
We recall some known facts about normal domains of attraction of stable distribution (see \cite{Feller}). Let $\xi$ belong to the normal domain of attraction of the stable distribution $S_\alpha(\sigma,\nu)$. Then, the following estimate holds for its characteristic function  as $\theta\to 0$ 
\begin{equation}\label{eq:fcstable2}
 \lambda(\theta)=\bbE\left[e^{i\theta\xi}\right]=\bar\lambda(\theta)+o(|\theta|^\alpha).
\end{equation}
where $\bar\lambda$ is given by (\ref{eq:fcstable1}).
Furthermore, in the case $0<\alpha<2$, the tail estimate (\ref{eq:heavy-tail}) implies that there exists $C>0$ such that for any $s>0$
\begin{equation}\label{eq:stablestim2}
{\rm Var}[\xi\ind_{\{|\xi|\leq s\}}]\leq Cs^{2-\alpha}\quad {\rm and} \quad  \bbE[|\xi|\ind_{\{|\xi|\leq s\}}]\leq Cs^{1-\alpha}.
\end{equation}

\subsection{Proof of Proposition \ref{prop:domain}}
The random measure defined by (\ref{eq:defmu}) is linked with the discretization of the space $\bbR^d$ by $h\bbZ^d$, $h>0$. 
Introduce the operator
$$
\psi_h :\left.\begin{array}{ccc} L^1_{loc} & \rightarrow & \bbR^{(\bbsZ^d)} \\ f & \mapsto &\left(h^{-d}\int_{h(k+I^d)}f(u)\d u\right)_{k\in\bbsZ^d} \end{array}\right..
$$ 
Define the random signed measure $\nu$ on $\bbsZ^d$  by 
\begin{equation}\label{eq:defnu}
\nu=\sum_{k\in\bbsZ^d}\xi_k \delta_k.
\end{equation}
It defines a random noise on the integrands set $\cF_\nu$ defined by
$$
\cF_\nu=\left\{f:\bbZ^d\to\bbR ;\ \\sum_{k\in\bbsZ^d}f_k\xi_k {\rm\ converges\ a.s.\ } \right\}
$$
where once again convergence means semi-convergence of the integral. 
We have the formal relation
$$
\mu_h[f]=\gamma_h h^{d}\nu[\psi_hf]
$$
from which we deduce that $f\in \cF_{\mu_h}$  if and only if $\psi_hf\in\cF_\nu$.
But we will see below that  $\cF_\nu=\ell^\alpha$ (see Lemma \ref{lem:Fnu}) and 
that $\psi_h(L^\alpha)\subset \ell^\alpha$ if $1\leq \alpha \leq 2$ and 
$\psi_h(\cD^\alpha)\subset \ell^\alpha$ if $0<\alpha < 1$ (see Lemma \ref{lem:lem1.2} in the Appendix).
This proves Proposition \ref{prop:domain}. \CQFD

\begin{lemme}\label{lem:Fnu}
Suppose $\alpha\in (0,2]$ and $\xi$ satisfies ${\bf H_\alpha}$. Then,  $F_\nu=\ell^\alpha$.
\end{lemme}
{\bf Proof of Lemma \ref{lem:Fnu}}\\
This is a direct application of Kolmogorov's three series Theorem (see \cite{Feller}).
The random series $\sum_{k\in\bbsZ^d} f_k\xi_k$ with independent summands converges almost surely if and only if for any $s>0$, the following
three numeric series converge: 
$$
\sum_{k\in \bbsZ^d} \bbP\left[|f_k\xi_k |>s\right]<\infty,\quad 
\sum_{k\in \bbsZ^d} {\rm Var}\left [f_k\xi_k \ind_{\{|f_k\xi_k |\leq s\}}\right]<\infty,
$$
and
$$
\sum_{k\in \bbsZ^d} \bbE\left[f_k\xi_k\ind_{\{|f_k\xi_k |\leq s\}}\right] \ \ {\rm converges}.
$$

Note that the convergence of the first series for small $s$ implies that $f_k\to 0$ as $k\to\infty$, 
so we do suppose $f_k\to 0$ as $k\to\infty$.
In the case $\alpha=2$, this implies
$${\rm Var}\left [f_k\xi_k \ind_{\{|f_k\xi_k |\leq s\}}\right]\sim_{k\to\infty} f_k^2\bbE[\xi_k^2]$$
and hence the second series converges if and only if $f\in \ell^2$. 
Suppose now $f\in \ell^2$. Cebycev's inequality implies
$$\bbP\left[|f_k\xi_k |>s\right]\leq f_k^2\bbE[\xi_k^2]s^{-2}$$ 
so that the first series converges. Applying Borel-Cantelli's Lemma, this in turn entails that 
the event $\limsup \{|f_k\xi_k |>s \}$ has probability $0$. Hence,  
the set $\{k;|f_k\xi_k |>s\}$ is almost surely finite, so that in the third series almost all terms 
vanish and the third series converges.

In the case $0<\alpha<2$, we deduce from equation (\ref{eq:heavy-tail}) that 
\begin{equation}\label{eq:stablestim1}
\bbP\left[|f_k\xi_k |>s\right]\sim (p+q)|f_k|^\alpha s^{-\alpha}
\end{equation}
and hence the first serie converges if and only if $f\in \ell^\alpha$.
Suppose now that $f\in\ell^\alpha$. Equation (\ref{eq:stablestim2}) implies the second and third series converge.
\CQFD

\subsection{Proof of Theorem \ref{theo:whitecase} and of Proposition \ref{prop:errorbound}}
{\bf Proof of Theorem \ref{theo:whitecase}}\\
Using Proposition \ref{prop:critere}, it is enough to prove the convergence of one dimensional distribution: for $f\in\cF^\alpha$, 
$$\mu_h[f] {\Rightarrow} W_{\alpha}[f]\quad \mbox{as}\quad h\to 0.$$ 
We will prove in fact a stronger result that will be useful in the sequel and we consider diagonal convergence.
Recall that $\cF^\alpha=L^\alpha$ if $1\leq \alpha \leq 2$ and $\cF^\alpha=\cD^\alpha$ if $0<\alpha<1$.
If $\alpha\in [1,2]$, $\cF^\alpha$ is a Banach space when endowed with the norm $|\!|.|\!|_{L^\alpha}$ and the notion of convergence $f_h\to f$ in $L^\alpha$ is clear. In the case $0<\alpha<1$, we need the following definition of convergence in $\cD^\alpha$.
\begin{df} 
We say that  $(f_h)_{h>0}$  converge to $f$ in $\cD^\alpha$ if the following two conditions hold:
\begin{itemize}
\item for any compact $K\subset \bbR^d$, $f_n{\bf 1}_K$ converges to $f{\bf 1}_K$ in $L^1$,
\item there is some $\eta>\alpha^{-1}d$ such that $f_h(x)=o(|x|^{-\eta} )$ and $f(x)=o(|x|^{-\eta})$ as $x\to\infty$ uniformly in $h>0$.
\end{itemize}
\end{df}

Theorem \ref{theo:whitecase} is then a direct consequence of the following Proposition.
\begin{prop}\label{prop:diagconv}Suppose $\xi$ satisfies ${\bf H_\alpha}$ and let $(f_h)_{h>0}$  converge to $f$ in $\cF^\alpha$ . Then, the following diagonal weak convergence holds:
$$\mu_h[f_h] {\Rightarrow} W_{\alpha}[f]\quad \mbox{as}\quad h\to 0.$$ 
\end{prop}

We prove convergence of the characteristic functions. Let $\theta\in\bbR$. From the definition of stable white noise, 
\begin{equation}\label{eq:fc-stableintegral}
\bbE\left[\exp(i\theta W_{\alpha}[f])\right]=\exp\left(-\sigma_f^\alpha|\theta|^\alpha(1-i\nu_f\varepsilon(\theta)\tan(\pi\alpha/2)) \right)
\end{equation}
with $\sigma_f$ and $\nu_f$ given by (\ref{eq:param}). On the other hand, we can rewrite 
\begin{equation}\label{eq:def-noise}
\mu_h[f_h]=\sigma^{-1}h^{\alpha^{-1}d} \sum_{k\in\bbsZ^d} \xi_k(\tilde\psi_h f_h)(hk),
\end{equation}
where $\tilde\psi_h :  L^1_{loc}\to L^1_{loc}$ is the linear functional defined by
\begin{equation}\label{eq:approx}
(\tilde\psi_h f)(x)=h^{-d}\int_{h([x]_h+I^d)}f(u)\d u,
\end{equation}
and $[x]_h=h[h^{-1}x]$. Note that $\tilde \psi_h f$ is the approximation of $f$ when the space $\bbR^d$ is discretized by $h\bbZ^d$ and the function $f$ is replaced by its mean value on each cell of the form $h(k+I^d)$, $k\in\bbZ^d$. 
With these notations, the characteristic function of $\mu_h[f_h]$ is given by
\begin{equation}\label{eq:fc-noise}
\bbE\left[\exp(i\theta\mu_h[f_h])\right]=\prod_{k\in\bbsZ^d} \lambda\left(\sigma^{-1}h^{\alpha^{-1}d}(\tilde\psi_h f_h)(hk)\theta \right).
\end{equation}
Here we have used equation (\ref{eq:def-noise}), the independence of the $\xi_k$'s and the almost sure convergence of the random serie (\ref{eq:def-noise}) implying the convergence of the above infinite product.\\ 

\noindent {\bf First step:} we begin to show that as $h\to 0$, 
\begin{equation}\label{eq:diff}
\bbE\left[\exp(i\theta\mu_h[f_h])\right]=\prod_{k\in\bbsZ^d} \bar\lambda\left(\sigma^{-1}h^{\alpha^{-1}d}(\tilde\psi_h f_h)(hk)\theta \right)+o(1).
\end{equation}
To see this, we estimate the difference and  use the following inequality : let $(z_i)_{i\in I}$ and $(z_i')_{i\in I}$ two families of complex numbers in $\bar{\bbD}$ such that the products $\prod_{i\in I}z_i$ and $\prod_{i\in I}z'_i$ are convergent, then
$$
\left|\prod_{i\in I}z'_i -\prod_{i\in I}z_i\right| \leq \sum_{i\in I} \left|z'_i - z_i\right|.
$$
Using equation (\ref{eq:fc-noise}), this yields 
\begin{eqnarray}
 & &\left|\ \bbE\left[\exp(i\theta\mu_h[f_h])\right]- \prod_{k\in\bbsZ} \bar \lambda \left( \sigma^{-1}h^{\alpha^{-1}d}(\tilde\psi_h f_h)(hk)\theta\right)\ \right| \nonumber \\
 &\leq& \  \sum_{k\in\bbsZ}  \ \left|\ \lambda \left(u_h(k)\right)- \bar\lambda \left(u_h(k) \right)\  \right|\label{eq:diff1}
\end{eqnarray}
with
$$
u_h(k)=\sigma^{-1}h^{\alpha^{-1}d}(\tilde\psi_h f_h)(hk)\theta.
$$
Equation (\ref{eq:fcstable2}) implies that the function $g$ defined by $g(0)=0$ and
$$
g(v)=|v|^{-\alpha}\left|\lambda(v)-\bar\lambda(v)\right|\ \ ,\ \ v\neq 0,
$$
is continuous and bounded so that for any $k\in\bbZ^d$,
$$
\left| \lambda \left( u_{h}(k)\right)-  \bar \lambda \left( u_{h}(k)\right) \right| = g(u_{h}(k))|u_{h}(k)|^\alpha .
$$
In order to obtain an uniform estimation, define the function $\tilde g:[0,+\infty)\rightarrow [0,+\infty)$ by
$$
\tilde g(u) =\sup_{|v|\leq u} |g(v)|.
$$
Note that $\tilde g$ is continuous, bounded and vanishes at $0$, and that for any $k\in \bbZ^d$ such that $|u_h(k)|\leq \varepsilon$,
\begin{equation}\label{eq:diff2}
\left| \lambda \left( u_{h}(k)\right)-  \bar \lambda \left( u_{h}(k)\right) \right|\leq \tilde g(\varepsilon)|u_{h}(k)|^\alpha .
\end{equation}
Let $\varepsilon>0$. Equations (\ref{eq:fc-noise}), (\ref{eq:diff1}) and (\ref{eq:diff2}) together yield
\begin{eqnarray*}
& &\left|\ \bbE\left[\exp(i\theta\mu_h[f_h])\right] - \prod_{k\in\bbsZ} \bar \lambda \left( \sigma^{-1}h^{\alpha^{-1}d}(\tilde\psi_h f_h)(hk)\theta\right)\ \right|\\
&\leq& \  \tilde g(\varepsilon)\sum_{k\in\bbsZ^d}\left|u_h(k) \right|^\alpha\ind_{\{|u_h(k)|\leq\varepsilon\}} + 2\sum_{k\in\bbsZ^d}\ind_{\{|u_h(k)|>\varepsilon\}}.
\end{eqnarray*}
Here we have used equation (\ref{eq:diff2}) to bound $|\bar\lambda(u_h(k))-\lambda(u_h(k))|$ from above whenever $|u_h(k)|\leq\varepsilon$, and the trivial upper bound $2$ otherwise.
Now we remark that 
$$
\sum_{k\in\bbsZ^d}   \left|u_h(k) \right|^\alpha \ind_{\{|u_h(k)|\leq\varepsilon\}}
\leq
\sum_{k\in\bbsZ^d}   \left|u_h(k) \right|^\alpha 
= \sigma^{-\alpha}|\theta|^\alpha |\!|\tilde\psi_h f_h|\!|_{L^\alpha}^\alpha 
$$
is bounded since from Corollary \ref{cor:cor1} $|\tilde\psi_h f_h|^\alpha \to  |f|^\alpha$ in $L^1$  whenever $f_h\to f$ in $\cF^\alpha$. By the continuity of $\tilde g$ in $0$,  $\tilde g(\varepsilon)$ is small for $\varepsilon$ small enough. Furthermore,
\begin{eqnarray*}
\sum_{k\in\bbsZ^d}\ind_{\{|u_h(k)|>\varepsilon\}}&=&h^{-d} \int_{\bbsR^d}\ind_{\{|\tilde\psi_hf_h(x)|>\varepsilon\sigma|\theta|^{-1}h^{-\alpha^{-1}d} \}}\d x\\
&\leq& \varepsilon^{-\alpha}\sigma^{-\alpha}|\theta|^\alpha |\!|(\tilde\psi_h f_h)\ind_{\{ |\tilde\psi_hf_h|>\varepsilon\sigma|\theta|^{-1}h^{-\alpha^{-1}d}\}}|\!|_{L^\alpha}^\alpha 
\end{eqnarray*}
and this quantity goes to $0$ as $h\to 0$ because Corollary \ref{cor:cor1} implies that $|\tilde\psi_h f_h|^\alpha \to  |f|^\alpha$ in $L^1$  whenever $f_h\to f$ in $\cF^\alpha$ and hence the familiy $|\tilde\psi_hf_h|^\alpha, h>0$ is uniformly integrable. 
These estimates imply equation (\ref{eq:diff}).\\

\noindent {\bf Second step:} we prove that as $h\to 0$,
\begin{equation}\label{eq:limit}
\prod_{k\in\bbsZ^d} \bar\lambda\left(\sigma^{-1}h^{\alpha^{-1}d}(\tilde\psi_h f_h)(hk)\theta \right)= \exp\left(-\sigma_f^\alpha|\theta|^\alpha(1-i\nu_f\varepsilon(\theta)\tan(\pi\alpha/2)) \right)+o(1).
\end{equation}
To see this, we take here advantage of the exponential form of $\bar\lambda$ and write the l.h.s. of (\ref{eq:limit}) as
\begin{eqnarray}
& &\exp\left( - |\theta|^\alpha h^{d}\sum_{k\in\bbsZ^d}|(\tilde\psi_h f_h)(hk)|^\alpha +i\nu \tan(\pi\alpha/2) \theta^{<\alpha>} h^{d}\sum_{k\in\bbsZ^d}((\tilde\psi_h f_h)(hk))^{<\alpha>}  \right) \nonumber\\
&=& \exp\left( - |\theta|^\alpha \int_{\bbsR^d} |(\tilde\psi_h f_h)(x)|^\alpha \d x +i\nu \tan(\pi\alpha/2) \theta^{<\alpha>} \int_{\bbsR^d} |(\tilde\psi_h f_h)(x)|^{<\alpha>} \d x  \right) \label{eq:calcul}
\end{eqnarray}
where $u^{<\alpha>}=\varepsilon(u)|u|^\alpha$ denotes the signed power function.

Corollary \ref{cor:cor1} implies that
$$
\int_{\bbsR^d} |(\tilde\psi_h f_h)(x)|^\alpha \d x=\int_{\bbsR^d} |f(x)|^\alpha \d x+o(1)
$$
and
$$
\int_{\bbsR^d} |(\tilde\psi_h f_h)(x)|^{<\alpha>} \d x=\int_{\bbsR^d} |f(x)|^{<\alpha >}\d x+o(1).
$$
This together with equation (\ref{eq:calcul}) implies equation (\ref{eq:limit}).\\

\noindent
{\bf Conclusion:} The convergence of the characteristic functions of $\mu_h[f_h]$ given by (\ref{eq:fc-noise}) to the characteristic function of $W[f]$ given by (\ref{eq:fc-stableintegral}) is a direct consequence of equations (\ref{eq:diff}) and (\ref{eq:limit}). This proves the weak convergence $\mu_h[f_h]\Rightarrow W[f]$. Proposition \ref{prop:diagconv} and Theorem \ref{theo:whitecase} are proved.\CQFD \\

\noindent
{\bf Proof of Corollary \ref{cor1}}\\
Note that formally $\tilde\mu_h[f]=\mu_h[f_h]$, where $f_h$ is the function defined by $f_h(x)=f([x]_h)$. The result follows from Proposition \ref{prop:diagconv} and from the fact that if $f\in\cF^\alpha\cap \cC(\bbR^d)$, then $f_h\to f$ in $\cF^\alpha$.\CQFD \\

\noindent
{\bf Proof of Proposition \ref{prop:errorbound}}\\
The proof is rather straightforward once we observe that 
$\bar\mu_h[f_M]-W_\alpha[f]=W_\alpha[\phi_hf_M -f]$ has the same distribution as $|\!|\phi_hf_M -f|\!|_{L_\alpha}X_\alpha$ where $X_\alpha$ has distribution $\cS_\alpha(1,\nu)$. Hence,
$$\bbE\left[ |\hat\mu_h[f_M]-W_\alpha[f]|^p\right]=  |\!|\tilde\psi_hf_M -f|\!|_{L_\alpha}^p\bbE[|X_\alpha|^p],$$ 
and $C_{\alpha,p}=\bbE[|X_\alpha|^p]$.\CQFD

\subsection{Proof of Theorem \ref{theo:poisson}}
We first prove that for any $\lambda>0$, the integrand set $\cF_{\mu_\lambda}$ contains $L^\alpha(E,\cE,m)$.
Let $f:E\to \bbR$ be such that $\int_{E} |f|^\alpha dm<\infty$. We have to prove that equation (\ref{eq:def_Fmulambda}) is satisfied and compute
\begin{eqnarray*} 
\int_{E\times\bbsR} |\xi f(e)|\wedge 1\ n_\lambda(d\xi,de)&=& \lambda \int_E \left(|f(e)|\bbE[|\xi|{\bf 1}_{|f(e)\xi| \leq 1}]+\bbP[|f(e)\xi|>1]\right) m(de) \\
&\leq& \lambda \int_E C |f(e)|^\alpha m(de) <\infty
\end{eqnarray*}
where the last inequality is a consequence of equations (\ref{eq:stablestim1}) and (\ref{eq:stablestim2}). This proves the required inclusion.

We now prove the convergence of the random measures to the stable random noise. Following Proposition \ref{prop:critere},
we compute the characteristic function of $N_\lambda[f]$ for $f\in L^\alpha(E,\cE,m)$:
$$\bbE\left[ \exp(i\theta N_\lambda[f])\right]=\exp\left(\int_{E\times\bbsR} \Psi(\theta\gamma_\lambda \xi f(e)) \lambda m(de)G(d\xi)\right), $$
 where $\Psi(u)=e^{iu}-1$. Performing integration with respect to $G$, the integral in the right hand side rewrites
$$\int_E \Psi_G\left(\theta\gamma_\lambda \xi f(e)) \right)\lambda m(de)G(d\xi),$$
with $\Psi_G(u)=\int_{\bbsR}\Psi(u\xi)G(d\xi)$. 
We now prove the following asymptotic result:
$$\lim_{\lambda\to\infty} \int_E \Psi_G\left(\theta\gamma_\lambda  f(e) \right)\lambda m(de)=-\int_E |\theta f(e)|^\alpha(1-\nu\varepsilon(\theta f(e))\tan(\pi\alpha/2)).$$
This is a consequence of Lebesgue's convergence Theorem. Indeed $\gamma_\lambda\to 0$ as $\lambda\to\infty$ and the asymptotic behaviour of $\Psi_G$ at $0$ given by equation (\ref{eq:fcstable2}) yields
$$\Psi_G(\theta\gamma_\lambda f(e))=-|\theta f(e)|^\alpha(1-\nu\varepsilon(\theta f(e))\tan(\pi\alpha/2))+o(1).$$
Furthermore, there is some $C>0$ such that $|\Psi_G(\theta)|\leq C|\theta|^\alpha$ for some $C>0$ and this implies the domination condition
$$ \lambda \Psi_G(\theta \gamma_\lambda f(e)) \leq C|\theta|^\alpha |f(e)|^\alpha$$
where the r.h.s. is integrable with respect to $m(de)$ and does not depend on $\lambda$.
Hence we have proved that 
$$\lim_{\lambda\to\infty}\bbE\left[ \exp(i\theta N_\lambda[f])\right]=\exp\left(-\int_{E} |\theta f(e)|^\alpha(1-\nu\varepsilon(\theta f(e))\tan(\pi\alpha/2)) m(de)\right), $$
and this is precisely the characteristic function of $W_\alpha[f]$. \CQFD

\subsection{Proof of Theorem \ref{theo:fraccase}}

In view of Lemma \ref{lem:Fnu}, the conditions $c\in\ell^\alpha$ and $\xi$ satisfies ${\bf H}_\alpha$ ensure  that $\hat\xi=\xi\ast c$ is well defined.
This convolution relation implies that the random measures $\mu_h$ and $\hat\mu_h$ defined by equations (\ref{eq:defmu}) and (\ref{eq:defhatmu}) respectively are linked by a simple relation. This is the object of the following Lemma. In order to unify the two cases considered in Theorem \ref{theo:fraccase}, we denote the integrands set by $\hat\cF_\alpha$ with $\hat\cF_\alpha=L^\alpha$ in the first case and $\hat\cF_\alpha=L^\alpha\cap L^1$ in the second case. We define $\check c\in\bbR^{\bbsZ^d}$ by $\check c_k=c_{-k},\ k\in\bbZ^d$.\\

\begin{lem}
Under the assumptions of Theorem \ref{theo:fraccase} (both cases), $\hat\cF^\alpha\subseteq \cF_{\hat\mu_h}$ and the relation  $\hat\mu_h[f]= \mu_h[\tilde{\psi}_h^c f]$ holds for  any $f\in\hat\cF^\alpha$ with
\begin{equation}\label{eq:OP3}
\tilde\psi_h^c f=\frac{\hat \gamma_h}{\gamma_h} \phi_h((\psi_h f)\ast \check c),
\end{equation}
where 
\begin{equation}\label{eq:OP2}
\psi_h :\left.\begin{array}{ccc} L^1_{loc} & \rightarrow & \bbR^{\bbsZ^d} \\ f & \mapsto &h^{-d}\left(\int_{h(k+I^d)}f(u)\d u\right)_{k\in\bbsZ^d} \end{array}\right.,
\end{equation}
and
\begin{equation}\label{eq:OP1}
\phi_h :\left.\begin{array}{ccc} \bbR^{\bbsZ^d}& \rightarrow &  L^1_{loc} \\ u=(u_k)_{k\in\bbsZ^d} & \mapsto & f:x\mapsto u_{[h^{-1}x]}\end{array}\right..
\end{equation}

\end{lem}

\noindent
{\bf Proof:}\\
We have formally,  
\begin{eqnarray*}
\hat \mu_h^{\xi}[f]&=&\hat\gamma_h \sum_{k\in\bbsZ^d} (\psi_hf)_k\hat\xi_k \\
&=& \hat\gamma_h \sum_{(k,l)\in\bbsZ^d\times \bbsZ^d } (\psi_hf)_k c_{k-l} \xi_l  \\
&=& \hat\gamma_h \sum_{l\in\bbsZ^d }  ((\psi_h f)\ast \check c)_l \xi_l\\
&=& \frac{\hat \gamma_h}{\gamma_h} \mu_h[\phi_h((\psi_h f)\ast \check c)].\\
\end{eqnarray*}
Furthermore, this formal computations are valid because the assumptions of Theorem \ref{theo:fraccase} entail $(\psi_h f)\ast \check c \in \ell^\alpha$ and hence $\tilde\psi_h^c f\in L^\alpha\subseteq \cF_\mu$. Recall indeed that Hölder inequality implies that $u\ast v\in \ell^\alpha$ as soon as $u\in \ell^1$ and $v\in\ell^\alpha$. Here in the first case (resp. in the second case) $c\in \ell^1$ and $\psi_h f\in \ell^\alpha$ (resp. $c\in \ell^\alpha$ and $\psi_h f\in \ell^1$). 
\CQFD
\ \\
\noindent
{\bf Proof of Theorem \ref{theo:fraccase}:}\\
 Theorem \ref{theo:fraccase} is a direct consequence of the above Lemma, of Theorem \ref{theo:whitecase} (or rather of Proposition \ref{prop:diagconv}  about diagonal convergence) and of the convergence of $\tilde\psi_h^c f$ in $L^\alpha$ as $h\to 0$. Lemma \ref{lem:lem1.3} in the appendix states indeed that $\tilde\psi_h^c f\to Cf$ converges to $Cf$ (resp. to $f\ast\check p$) in the first case (resp. second case). \CQFD

\section*{Appendix: some results on space functions}


The following Lemma gathers some useful properties of the linear operators $\tilde\psi_h,\psi_h$ and $\phi_h$ defined by equations (\ref{eq:approx}), (\ref{eq:OP2}) and (\ref{eq:OP1}) respectively. It is standard material for $\alpha\geq 1$ but some extra care is needed when $0<\alpha<1$.

\begin{lemme}\label{lem:lem1.2}
\begin{enumerate}
\item The relations $\tilde\psi_h=\phi_h\circ \psi_h$, $\psi_h\circ \phi_h=Id_{\bbsR^{\bbsZ^d}}$ and $\tilde\psi_h\circ \tilde\psi_h =\tilde\psi_h $ hold.
\item $\phi_h(\ell^\alpha)\subset L^\alpha$ and for any $u\in \ell^\alpha$, $|\!|\phi_h u|\!|_{L^\alpha}=h^{d/\alpha}|\!|u|\!|_{\ell^\alpha}$.
\item If $\alpha\geq 1$, $\psi_h(L^\alpha)=\ell^\alpha$ and  $\tilde\psi_h(L^\alpha)\subset L^\alpha$ and $\tilde\psi_h$ induces a linear projection on $L^\alpha$. 
Furthermore, if $f_h\to f$ in $L^\alpha$ as $h\to 0$, then $\tilde\psi_h(f_h)\to f$ in $L^\alpha$.
\item If $0<\alpha<1$,  $\psi_h(\cD_\alpha)\subset \ell^\alpha$, 
  and $\tilde\psi_h(\cD_\alpha)\subset (\cD_\alpha)$. 
Furthermore, if $f_h\to f$ strongly in $\cD_\alpha$ as $h\to 0$, then $\tilde\psi_h(f_h)\to f$   strongly in $\cD_\alpha$.
\end{enumerate}
\end{lemme} 
\ \\
\noindent
{\bf Remarks:} Extra care is needed in the case $0<\alpha<1$ because  $\tilde\psi_h(L^1_{loc}\cap L^\alpha)$ is not included in $L^\alpha$ as the following example shows: for $\gamma>0$ and $\delta\in\bbR$  consider the function $f\in L^1_{loc}$ defined by $f(x)=k^\delta $ if $x\in [k, k+k^{-\gamma}]$ for  $k\geq 1$, and $f(x)=0$ otherwise. Then it is easily seen that $f\in L^\alpha$ if and only if $\alpha\delta-\gamma<-1$, and that $\tilde\psi_1 f\in L^\alpha$ if and only if $\alpha(\delta-\gamma)<-1$. Comparing these to conditions, we can even find $f\in L^\alpha\cap L^1$ such that $\tilde\psi_1 f\notin L^\alpha$. That's why we need some stronger condition on $f$ in the case $0<\alpha<1$ and we introduce the space of rapidly decaying functions $\cD_\eta$.

\ \\

\noindent
{\bf Proof of Lemma \ref{lem:lem1.2}:}\\
\begin{enumerate}
\item This is routine verification using the definitions.
\item Let $u\in\ell^\alpha$. Then,
\begin{eqnarray*}
\int_{bbsR^d}|\phi_h u(x)|^\alpha dx&=& \int_{bbsR^d}|u_{[h^{-1}x]}|^\alpha dx\\
&=& h^{d}\sum_{k\in\bbsZ^d} |u_k|^\alpha, 
\end{eqnarray*}
and this implies $\phi_hu\in L^\alpha$ and $|\!|\phi_h u|\!|_{L^\alpha}=h^{d/\alpha}|\!|u|\!|_{\ell^\alpha}$.
\item Let $\alpha\geq 1$ and $f\in L^\alpha$. Using Jensen inequality with the convex function $y\mapsto |y|^\alpha$, 
\begin{eqnarray*}
|\!|\psi_h f|\!|_{\ell^\alpha}^\alpha&=&h^{-d}\sum_{k\in\bbsZ^d}\left|h^{-d}\int_{h(k+I^d)} f(u)du \right|^\alpha\\
&\leq& \sum_{k\in\bbsZ^d} h^{-d}\int_{h(k+I^d)} |f(u)|^\alpha du\\
&=& h^{-d}|\!|f|\!|_{L^\alpha}^\alpha.
\end{eqnarray*}
This implies that $\psi_h f\in \ell^\alpha$ and $|\!|\psi_h f|\!|_{\ell^\alpha}\leq h^{-d/\alpha}|\!|f|\!|_{L^\alpha}$.
Together with point 1. and 2., we easily prove that $\tilde\psi_h$ induces a linear projection on $L^\alpha$. \\
We now prove that $\tilde\psi_h$ converge weakly as a linear operator to the idendity.
First consider the case when $f$ is continuous with compact support. Then from the continuity of $f$, for every $x\in\bbR\d$ $(\tilde\psi_hf)(x)\rightarrow f(x)$ as $h\to 0$. The assumption that $f$ is compactly supported implies that Lebesgue's convergence theorem applies and $\psi_h f\to f$ in $L^\alpha$.\\
The space of continuous compactly supported function is dense in $L^\alpha$ and the $\psi_h$ are (uniformly) bounded linear operator. Hence the convergence $\psi_h f\to f$ holds for every $f\in L^\alpha$. Diagonal convergence is a consequence of the equicontinuity of the $\psi_h$ which all have linear norm less than $1$.
\item Let $f\in \cD^\alpha$. From the definition of $\cD^\alpha$, there is some integer $M>0$ and $C>0$ and some 
$\eta>\alpha^{-1}d$ such that if $|x|_\infty>M$ then $|f(x)|\leq C|x|_\infty^{-\eta}$. As a consequence, if $h\leq 1$ and $|x|_\infty\geq M+h$, we have $h([h^{-1}x]+I^d)\subseteq B(O,M)^c$ and 
$$|\tilde\psi_h f(x)|^\alpha=\left|h^{-d}\int_{h([h^{-1}x]+I^d)}f(u)\d u \right|^\alpha\leq h^{(1-\alpha)d}(|x|_\infty-h)^{-\alpha\eta}.$$
This implies $\int_{\bbsR^d}|\tilde\psi_hf(x)|^\alpha\d x<\infty$ since $\alpha\eta>d$ and hence $\tilde\psi_h f\in L^\alpha$.

Now let $f_h\rightarrow f$ strongly in $\cD^\alpha$. From the above discussion and the definition of strong convergence in $\cD^\alpha$, we  see that a similar inequality holds uniformly in $h\leq 1$: there is some  $M>0$,  $C>0$ and $\eta>\alpha^{-1}d$ such that if $|x|_\infty>M$ then $|f(x)|\leq C|x|_\infty^{-\eta}$  and also for any $h\leq 1$
$|f_h(x)|\leq C|x|_\infty^{-\eta}$. 
As a consequence, for any $\epsilon>0$, there is some $M>0$ such that for any $h\leq 1$
\begin{equation}\label{eq:maj-unif}
\int_{|x|_\infty>M} |\tilde\psi_hf_h(x)|^\alpha \d x  \leq \epsilon
\end{equation}
On the other hand, for fixed $M>0$, we  prove that $|\tilde\psi_hf_h|^\alpha {\bf 1}_{B(0,M)}$ is uniformly integrable.
Let $K=B(0,M+1)$. Then $f_h{\bf 1}_K\to f{\bf 1}_K$ in $L^1$, and  point 3. implies that $\tilde\psi_h (f_h {\bf 1}_K) \to f{\bf 1}_K$ in $L^1$. Since the restrictions of $ \tilde\psi_h (f_h {\bf 1}_K)$ and $\tilde\psi_h f_h$ to $B(0,M)$ are equal, this implies $(\psi_h  f_h) {\bf 1}_{B(0,M)}$ converge to  $f{\bf 1}_{B(0,M)}$ in $L^1$. 
This in turn implies that $|\tilde\psi_h  f_h|^\alpha {\bf 1}_{B(0,M)}$ is bounded in  $L^p$ for $p=\alpha^{-1}>1$, and hence uniformly integrable. This together with equation (\ref{eq:maj-unif}) implies the uniform integrability of the sequence $|\tilde\psi_nf_n|^\alpha$.
\CQFD
\end{enumerate}

\begin{corol}\label{cor:cor1}
Let $\alpha \in (0,2]$. Suppose $(f_h)_{h>0}$ converge to $f$ in $\cF^\alpha$ as $h\to 0$. Then $|\tilde\psi_h(f_h)|^\alpha \to |f|^\alpha$ in $L^1$ as $h\to 0$ and also $(\tilde\psi_h(f_h))^{<\alpha>} \to f^{<\alpha>}$ in $L^1$
\end{corol}

\noindent
{\bf Proof of Corollary \ref{cor:cor1}:}\\
We consider the convergence of $(\tilde\psi_hf_h(x))^{<\alpha>}$ to $f^{<\alpha>}$ in $L^1$, the other case being treated in the same way. First of all, as above, for any $\epsilon>0$ there is some $M>0$ such that for any $h\leq 1$
 \begin{equation}\label{eq:maj-unif2}
\int_{|x|_\infty>M} \left|(\tilde\psi_hf_h(x))^{<\alpha>}-f(x)^{<\alpha>}\right| \d x  \leq \epsilon.
\end{equation}
For fixed $M>0$, the application $L^1(B(0,M))\to \bbR, f\mapsto \int_{B(0,M)} f(x)^{<\alpha>}\d x$ is continuous.
Hence the convergence  $(\psi_h  f_h) {\bf 1}_{B(0,M)}\to f{\bf 1}_{B(0,M)}$ in $L^1$ implies
$$\lim_{n\to\infty} \int_{B(0,M)} (\psi_h  f_h)(x)^{<\alpha>}\d x= \int_{B(0,M)} (\psi_h  f_h)(x)^{<\alpha>}\d x.$$
Together with equation (\ref{eq:maj-unif2}), this yields the result.
\CQFD

\ \\

We now consider the convergence of $\tilde\psi_h^c f$ defined by (\ref{eq:OP3}).
\noindent
\begin{lem}\label{lem:lem1.3} Under the assumptions of Theorem \ref{theo:fraccase}
\begin{enumerate} 
\item In the first case, for any $f\in L^\alpha$, $\tilde\psi_h^c f\to Cf$ in $L^\alpha$ as $h\to 0$.
\item In the second case, for any $f\in L^1\cap L^\alpha$, $\tilde\psi_h^c f\to f\ast\check p$ in $L^\alpha$ as $h\to 0$.
\end{enumerate}
\end{lem}

\noindent
{\bf Proof of Lemma \ref{lem:lem1.3}:}\\
\begin{enumerate}
\item Consider first the case when $c\in\ell^1$. We easily see that
$$\tilde\psi_h^c f=\sum_{k\in\bbsZ^d}c_{k} \tilde\psi_hf(\cdot+ hk ).$$
For each fixed $k$, the sequence of function $\tilde\psi_hf(\cdot+ hk )$ converge to $f$ in $L^\alpha$ (use the fact that according to Lemma \ref{lem:lem1.2} point 3, $\tilde\psi_hf\to f$ in $L^\alpha$ and that the translation operator $f\to f(.+hk)$ is continuous and converge to identity as $h\to 0$). Then for any fixed $M$, 
$$\sum_{|k|\leq M}c_{k} \tilde\psi_hf(\cdot+ hk )\to (\sum_{|k|\leq M}c_{k}) f $$
in $L^\alpha$ as $h\to 0$.  Furthermore, since $c\in\ell^1$, the remainder 
$$|\!|\sum_{|k|> M}c_{k} \tilde\psi_hf(\cdot+ hk )|\!|_{L^\alpha} \leq \left(\sum_{|k|> M}|c_{k}|\right)  |\!|f|\!|_{L^\alpha} $$
is small when $M$ is large. Thanks to these estimates, we prove that $\tilde\psi_h^c f\to Cf$ in $L^\alpha$.
\item Consider the second case when $c$ is such that
$$\lim_{t\to\infty} \sup_{|x|=1} \left|t^\beta c_{[tx]}-p(x)\right|=0,$$
and note that $\tilde\psi_h^c f$ can be rewritten as 
\begin{equation}\label{eq:quasiconv}
\tilde\psi_h^c f(x)= \int_{\bbsR^d}\tilde\psi_h f([x]_h-y) \check p_h(y)  \d y=((\tilde\psi_hf) \ast \check p_h)([x]_h)
\end{equation}
with
$$\check p_h(y)=h^{-d}\frac{\hat\gamma_h}{\gamma_h}c_{-[h^{-1}y]}=h^{-\beta}\check c_{[h^{-1}y]}.  $$
Using equation \ref{eq:asympt}, we see that for any $y\neq 0$, $\check p_h(y)\to \check p(y)$ as $h\to 0$.\\
More precisely, equation \ref{eq:asympt} entails that $\check p_h1_{B}\to \check p_h1_{B}$ in $L^1$ and $\check p_h1_{B^c}\to \check p_h1_{B^c}$ in $L^\alpha$, where $B$ denotes the unit ball $B(0,1)$ in $\bbR^d$ and $B^c$ its complementary set. Using Lemma \ref{lem:lem1.2} and the assumption $f\in L^1\cap L^\alpha$, we have also $\tilde\psi_h f\to f$ in $L^1$ and in $L^\alpha$. The bilinear application $L^1\times L^\alpha \to L^\alpha , (g_1,g_2)\mapsto g_1\ast g_2 $ is continuous and
$$\left|\!\left| g_1\ast g_2\right|\!\right|_{L^\alpha} \leq \left|\!\left| g_1 \right|\!\right|_{L^1}\left|\!\left|  g_2\right|\!\right|_{L^\alpha} .$$ This entails $(\tilde\psi_h f )\ast (\check p_h1_{B}) \to f\ast (\check p1_{B})$ in $L^\alpha$ and $(\tilde\psi_h f) \ast (\check p_h1_{B^c}) \to f\ast (\check p1_{B^c})$, and finally $(\tilde\psi_h f) \ast \check p_h \to f\ast \check p$ in $L^\alpha$.

At last, we prove that $\tilde\psi_h^c f\to f\ast\check p$ in $L^\alpha$. We have indeed
\begin{eqnarray*}
& &\left|\!\left| \tilde\psi_h^c f -(\tilde\psi_h f)\ast \check p_h \right|\!\right|_{L^\alpha}^\alpha \\
&\leq & \int_{\bbsR^d} \left|  \int_{\bbsR^d}|\tilde\psi_h f([x]_h-y)-\tilde\psi_h f(x-y)|  |\check p_h(y)|  \d y \right|^\alpha \d x \\
&\leq &|\!|\check p_h 1_B|\!|_{L^1}^{\alpha}  \sup_{0\leq u\leq h} |\!|\tilde\psi_h f(.+u)-\tilde\psi_h f|\!|_{L^\alpha}^{\alpha}+  |\!|\check p_h 1_{B^c}|\!|_{L^\alpha}^{\alpha}  \sup_{0\leq u\leq h} |\!|\tilde\psi_h f(.+u)-\tilde\psi_h f|\!|_{L^1}^{\alpha} 
\end{eqnarray*}
and these quantities vanish as $h\to 0$.
\CQFD
\end{enumerate}

\section*{Acknowledgements}
I am grateful to Professor Serge Cohen for fruitful discussions at the origin of the limit theorems for the Lévy Brownian motion on the sphere and related processes.


\begin{thebibliography}{2}

\bibitem{BMS} {\sc Biermé H.}, {\sc Meerschaert M.M.} and {\sc Scheffler H.-P.} {\it Operator scaling stable random fields}. Preprint available on arXiv at http://arxiv.org/abs/math/0602664v1.

\bibitem{CM89} {\sc Cambanis S.} and {\sc Maejima M.} (1989) {\it Two classes of self-similar stable processes with stationary increments.} {\rm Stochastic Processes and their Applications} 32:305--329.

\bibitem{CLL08} {\sc Cohen S.}, {\sc Lacaux C.} and {\sc Ledoux M.} (2008) {\it A general framework
for simulation of fractional fields}. {\rm Stochastic Process. Appl.}, 118(9):1489--1517.

\bibitem{CS} {\sc Cohen, S.} and {\sc Samorodnitsky, G.} (2006) {\rm Random rewards, fractional
Brownian local times and stable self-similar processes.} {\it Ann. Appl.
Probab.,} 16(3):1432--1461.

\bibitem{DOT03}{\sc Doukhan P.}, {\sc Oppenheim G.} and {\sc Taqqu M.S.} editors  (2003) {\it Theory and Applications of Long-Range Dependence}.  ISBN 0-8176-4168-8. {\rm Birkhäuser, Boston}. 

\bibitem{Feller} {\sc Feller W.} (1966) {\it An Introduction to Probability Theory and its Applications.} Vol. 2. 
{\rm Wiley}. 

\bibitem{K86}{\sc Kallenberg O.} (1986) {\rm Random Measures}, 4th ed. (187 pp). {\rm Akademie-Verlag and Academic Press, Berlin and London}.

\bibitem{KT94} {\sc Kokoszka P.} and {\sc Taqqu M.S.} (1994) {\rm Infinite variance stable ARMA processes}  {\it Journal of Time Series Analysis},  15:203--220.

\bibitem{KT94b}{\sc Kokoszka P.} and {\sc Taqqu M.S.} (1994) {\rm New classes of self-similar symmetric stable random fields}. {\it Journal of Theoretical Probability}, 7:527--549.

\bibitem{KT95}{\sc Kokoszka P.} and {\sc Taqqu M.S.} (1995) {\rm Fractional ARIMA with stable innovations}. {\it Stochastic Processes and their Applications}, 60:19--47.

\bibitem{KT96}{\sc Kokoszka P.} and {\sc Taqqu M.S.} (1996) {\rm Infinite variance stable moving averages with long memory}. {\it Journal of Econometrics}, 73:79--99.
 
\bibitem{KTU81} {\sc Kubo I.}, {\sc Takenaka S.} and {\sc Urakawa H.} (1981) {\it Brownian motion parametrized with metric space of constant curvature}. {\rm Nagoya Math. J.}  82:131-140. 

\bibitem{Ma83} {\sc Maejima M. } (1983) {\it On a class of self-similar processes}. {\rm Zeitschrift für Wahrscheinlichkeitstheorie und verwandte Gebiete.} 62:235--245.

\bibitem{Marouby08} {\sc Marouby, M.} {\it Simulation of local time stable motion}. Preprint available on arXiv at http://arxiv.org/abs/0712.3210v2.

\bibitem{Ro87} {\sc Rosinski J.} (1989) {\it On Path Properties of Certain Infinitely Divisible
Processes} {\rm Stochastic Processes and their Applications} 33(1):73--87. 

\bibitem{Ro90} {\sc Rosinski J.} {\it On Series Representations of Infinitely Divisible Random
Vectors}. {\rm The Annals of Probability}, 18(1):405--430.

\bibitem{Ro01} {\sc Rosinski J.} (2001) {\it Series representations of Lévy processes from the
perspective of point processes}. In Lévy processes, pages 401--415. {\rm Birkhauser}, Boston, MA. 

\bibitem{ST94} {\sc Samorodnitsky G.} and {\sc Taqqu M.S.} (1994) {\it Stable Non-Gaussian Random Processes: Stochastic Models with Infinite Variance.}, ISBN 0-412-05171-0, {\rm Chapman and Hall, New York}. 

\bibitem{TW83} {\sc Taqqu M.S.} and {\sc Wolpert R. } (1983) {\it Infinite variance self similar processes subordinate to a Poisson measure}. {\rm Zeitschrift für Wahrscheinlichkeitstheorie und verwandte Gebiete.} 62:53--72.



\end{thebibliography}
\end{document}